\documentclass[12pt]{article}

\usepackage{algorithm}				
\usepackage[noend]{algpseudocode}	
\usepackage{amsmath}				
\usepackage{amssymb}				
\usepackage{amsthm}					
\usepackage{caption}
\usepackage{color}
\usepackage[margin=1in]{geometry}
\usepackage[hidelinks]{hyperref}
\usepackage{titling}

\usepackage{mathtools}				
\usepackage{tikz}					
\usetikzlibrary{calc, intersections}
\newcommand{\drawrightangle}[4][.3cm]{
	\coordinate (tempa) at ($(#3)!#1!(#2)$);
	\coordinate (tempb) at ($(#3)!#1!(#4)$);
	\coordinate (tempc) at ($(tempa)!0.5!(tempb)$);
	\draw (tempa) -- ($(#3)!2!(tempc)$) -- (tempb);
}

\theoremstyle{plain}				
\newtheorem{theorem}{Theorem}
\newtheorem{lemma}[theorem]{Lemma}

\theoremstyle{definition}			
\newtheorem{definition}[theorem]{Definition}
\newtheorem{example}[theorem]{Example}

\theoremstyle{remark}				
\newtheorem*{remark}{Remark}

\newcommand{\divides}{\mid}
\newcommand{\from}{\colon}
\newcommand{\seqnum}[1]{\href{http://oeis.org/#1}{#1}}

\thanksmarkseries{alph}

\title{Primitive Two-Dimensional Words and Iterated Pedal Triangles via Symbolic Coding}
\author{
	Taylor J. Smith \thanks{Department of Computer Science, St.\ Francis Xavier University, Antigonish, Nova Scotia, Canada. Email: \href{mailto:tjsmith@stfx.ca}{\texttt{tjsmith@stfx.ca}}.}
}
\date{\today}

\begin{document}


\maketitle

\begin{abstract}
The notion of a two-dimensional word arises naturally in the study of combinatorics on words, while the iterative construction of pedal triangles results in a rich dynamical system in the study of geometry. At first, these two classes of objects seem to be unrelated. However, it is known that for all $n \geq 1$, the number of primitive two-dimensional words of dimension $2 \times n$ over a binary alphabet agrees with the number of triangles whose first similar pedal triangle is their $n$th pedal triangle. We construct a finite four-symbol coding of the sorted pedal map and use the resulting branch itineraries to give a bijection between these two classes.

\medskip


\medskip

\noindent\textit{MSC2020 classes:} 37B10 (primary); 51M15, 68Q45 (secondary).
\end{abstract}


\section{Introduction}\label{sec:introduction}

A \emph{two-dimensional word} is an array or matrix of symbols drawn from some finite alphabet $\Sigma$, and we say that a two-dimensional word with $m$ rows and $n$ columns has \emph{dimension} $m \times n$. If we cannot decompose a given two-dimensional word $W$ into some repeated concatenation of a smaller two-dimensional word, then we say that $W$ is a \emph{primitive} two-dimensional word. There exists a formula to enumerate all primitive two-dimensional words of dimension $m \times n$ over a binary alphabet $\Sigma$; call this value $\psi_{2}(m,n)$~\cite{GamardRichommeShallitSmith2017Periodicity, Smith2017MMathThesis}.

A \emph{pedal triangle} is constructed by starting with a triangle $T$, taking the three points where the altitudes dropped at each vertex of $T$ intersect with the opposite (possibly extended) side, and joining these points. We can iterate the process of constructing pedal triangles in the usual way: we construct the $n$th pedal triangle of $T$ by constructing the pedal triangle of the $(n-1)$st pedal triangle of $T$; the `zeroth' pedal triangle is $T$ itself. If $T$ is similar to its $n$th pedal triangle, where triangle similarity is determined by internal angles, then we say that $T$ has \emph{pedal period dividing $n$}. If $n$ is the least natural number with this property, then we say $T$ has \emph{exact pedal period $n$}. There exists a formula to enumerate all triangles having exact pedal period $n$; call this value $\chi(n)$~\cite{Valyi1903Fusspunktdreiecke}.

These two objects---one from formal language theory and one from geometry---seem completely unrelated, yet the author's master's thesis~\cite[p.\ 33]{Smith2017MMathThesis} observes that the values $\psi_{2}(2,n)$ and $\chi(n)$ coincide for all $n \geq 1$ and states that ``it is not known whether there exists a connection [\,\dots] but the fact that the two enumerations agree suggests some sort of bijection between these objects." In this paper, we use ideas from symbolic dynamics to establish such a bijection. Namely, we construct a finite four-symbol coding of the sorted pedal map and show that the branch itineraries of triangles having exact pedal period $n$ correspond to primitive words of length $n$ avoiding a certain form. We then apply a substitution to convert branch itineraries into primitive two-dimensional words of dimension $2 \times n$.


\section{Preliminaries}\label{sec:preliminaries}

In this section, we fix some notation and establish key properties about the objects studied in this paper. For more details about two-dimensional formal language theory in particular, see the article by Giammarresi and Restivo~\cite{GiammarresiRestivo19972DLanguages} or the book by Rosenfeld~\cite{Rosenfeld1979PictureLanguages}.


\subsection{Two-Dimensional Words}\label{subsec:2Dwords}

Let $\Sigma = \{0, 1, \dots, k-1\}$ denote a finite $k$-ary \emph{alphabet} for $k \geq 1$. In one dimension, a word $w = a_{0}a_{1} \dots a_{n-1}$ is a sequence of concatenated symbols $a_{i} \in \Sigma$, and the \emph{length} of $w$ is denoted $|w| = n$. A \emph{language} over $\Sigma$ is then a set of words over $\Sigma$.

Just as we may concatenate symbols to form words, we may concatenate words to form longer words. If we concatenate a word $w$ with itself $m$ times, we obtain the \emph{word power} $w^{m}$. A word $w$ is said to be \emph{periodic} if it can be written as $w = v^{m}$ where $v$ is a nonempty word and $m \geq 2$. If $w$ is not periodic, then it is said to be \emph{primitive}.

In two dimensions, we may develop analogous definitions of words, lengths of words, and languages. For $m \geq 1$ and $n \geq 1$, a \emph{two-dimensional word}
\begin{equation*}
W = 
\begin{bmatrix}
a_{0,0}		& a_{0,1}	& \dots		& a_{0,n-1} \\
a_{1,0}		& a_{1,1}	&			& a_{1,n-1} \\
\vdots		& 			& \ddots	& \vdots \\
a_{m-1,0}	& a_{m-1,1}	& \dots		& a_{m-1,n-1} \\
\end{bmatrix}
\end{equation*}
is a map from $\{0, 1, \dots, m-1\} \times \{0, 1, \dots, n-1\}$ to $\Sigma$, and we say that $W$ has \emph{dimension} $m \times n$. In this way, a two-dimensional word is similar to an array or matrix of symbols over $\Sigma$. The \emph{size} (i.e., `length') of the two-dimensional word $W$ is denoted $|W| = mn$, and a \emph{two-dimensional language} over $\Sigma$ is likewise a set of two-dimensional words over $\Sigma$. Note that a two-dimensional word of dimension $1 \times n$ for any $n \geq 1$ is just a usual (one-dimensional) word.

We denote the special language of all two-dimensional words over $\Sigma$ having dimension $m \times n$ by the notation $\Sigma^{m \times n}$. Furthermore, given a two-dimensional word $W$, we may index into $W$ to retrieve the symbol at position $(i,j)$ by writing $W[i,j]$, where $0 \leq i < m$ and $0 \leq j < n$.

Just as we may concatenate words in one dimension in the usual manner, we may concatenate two-dimensional words, but we must take care both to distinguish in which dimension we perform the concatenation and to match dimensions appropriately. Suppose we have a pair of two-dimensional words $W \in \Sigma^{m_{1} \times n_{1}}$ and $V \in \Sigma^{m_{2} \times n_{2}}$. If $n_{1} = n_{2} = n$, then the \emph{row concatenation} of $W$ and $V$ is the two-dimensional word of dimension $(m_{1} + m_{2}) \times n$ produced by joining $W$ and $V$ in such a way that the topmost row of $V$ is adjacent to the bottommost row of $W$. On the other hand, if $m_{1} = m_{2} = m$, then the \emph{column concatenation} of $W$ and $V$ is the two-dimensional word of dimension $m \times (n_{1} + n_{2})$ produced by joining $W$ and $V$ in such a way that the leftmost column of $V$ is adjacent to the rightmost column of $W$.

By repeatedly concatenating two-dimensional words, we may generalize the notion of a word power to two dimensions. Suppose $W \in \Sigma^{m \times n}$. Then the \emph{$p \times q$ power} of $W$, denoted $W^{p \times q}$, is the two-dimensional word of dimension $pm \times qn$ with the property that $W^{p \times q}[i,j] = W[i \bmod m, j \bmod n]$.

Lastly, defining two-dimensional word powers allows us to define primitivity and periodicity in two dimensions. A two-dimensional word $W$ is said to be \emph{periodic} if it can be written as $W = V^{p \times q}$ where $V$ is a nonempty two-dimensional word and $p \geq 2$ or $q \geq 2$ (or both). If $W$ is not periodic, then it is said to be \emph{primitive}.


\subsection{Pedal Triangles}\label{subsec:pedaltriangles}

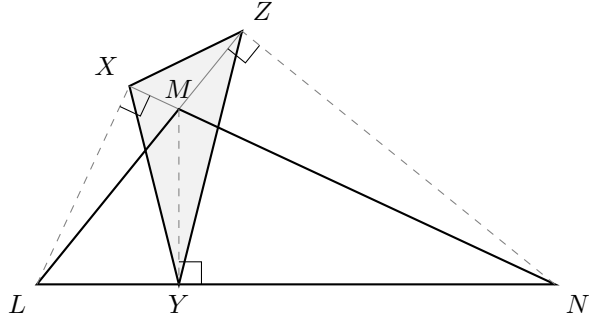
\begin{figure}[t]
\centering
\begin{tikzpicture}
	\def\A{3}
	\def\angA{25} 
	\def\angB{51} 
	\pgfmathsetmacro{\angC}{180-\angA-\angB}
	\pgfmathsetmacro{\d}{\A/sin(\angA)}
	\pgfmathsetmacro{\C}{\d*sin(\angC)}
	\coordinate (P) at (0,0);
	\coordinate (Q) at (\angB:\A);
	\coordinate (R) at (0:\C);
	\draw[black, thick, fill=gray!10] ($(P)!(Q)!(R)$) -- ($(Q)!(R)!(P)$) -- ($(R)!(P)!(Q)$) -- cycle;
	\draw[gray] (P) -- ($(Q)!(R)!(P)$);
	\draw[gray] ($(R)!(P)!(Q)$) -- (Q);
	\draw[black, thick] (P) -- (Q) -- (R) -- cycle;
	\draw[gray, dashed] (P) -- ($(R)!(P)!(Q)$);
	\draw[gray, dashed] (Q) -- ($(P)!(Q)!(R)$);
	\draw[gray, dashed] (R) -- ($(Q)!(R)!(P)$);
	\draw (P) node[anchor=north east]{\footnotesize $L$};
	\draw (Q) node[anchor=south]{\footnotesize $M$};
	\draw (R) node[anchor=north west]{\footnotesize $N$};
	\draw ($(R)!(P)!(Q)$) node[anchor=south east]{\footnotesize $X$};
	\draw ($(P)!(Q)!(R)$) node[anchor=north]{\footnotesize $Y$};
	\draw ($(Q)!(R)!(P)$) node[anchor=south west]{\footnotesize $Z$};
	\drawrightangle{P}{$(R)!(P)!(Q)$}{R}
	\drawrightangle{Q}{$(P)!(Q)!(R)$}{R}
	\drawrightangle{R}{$(Q)!(R)!(P)$}{P}
\end{tikzpicture}
\caption{A pedal triangle (in grey) constructed from a triangle (in white)}
\label{fig:pedaltriangle}
\end{figure}

Let $T$ be a triangle formed by joining vertices $L$, $M$, and $N$. If we drop the altitude from vertex $L$ to intersect the (possibly extended) line segment $MN$ and mark this point of intersection as $X$, then do the same for vertices $M$ and $N$ intersecting line segments $LN$ and $LM$ to mark points $Y$ and $Z$, respectively, then the triangle formed by joining $X$, $Y$, and $Z$ is known as the \emph{pedal triangle} of $T$. Figure~\ref{fig:pedaltriangle} depicts an example construction of a pedal triangle.

\begin{remark}
The geometric object produced according to our definition of a pedal triangle ought to more properly be referred to as an \emph{orthic triangle}. The most general definition of a pedal triangle is one constructed by joining the three points where the perpendiculars drawn from an arbitrary point $P$ intersect each side of a triangle; our definition arises by taking $P$ to be the orthocenter of the triangle. This slight abuse of terminology is unfortunate, but consistent with much of the literature on pedal triangles.
\end{remark}

We may iterate the process of constructing pedal triangles in the following way. If we take $T_{0} = T$ to be the original triangle and denote its pedal triangle by $T_{1} = T_{0}^{\prime}$, then in general we may construct and denote the $n$th pedal triangle of $T$ by taking $T_{n} = T_{n-1}^{\prime}$ for any $n \geq 1$.

Note that if we construct the pedal triangle of a right-angled triangle, then two of the vertices of the pedal triangle will coincide and the construction will produce a degenerate triangle (i.e., a line segment). We therefore refer to right-angled triangles, as well as any triangle for which some pedal iterate is right-angled, as being \emph{pedally degenerate}.

If a triangle $T$ is similar to its $n$th pedal triangle $T_{n}$, then we say that $T$ has \emph{pedal period dividing $n$}. Moreover, if $T$ is similar to $T_{n}$ but not to any other pedal triangle $T_{m}$ where $1 \leq m < n$, then we say that $T$ has \emph{exact pedal period $n$}.

Our objects of interest in this paper are triangles whose pedal iterates are similar to the original triangle. Since triangle similarity is determined by interior angles, we will represent triangles as triples of angles. Suppose a triangle $T$ has angles $A$, $B$, and $C$, where $A + B + C = \pi$. For the sake of convenience, we will normalize these angles by writing $a = A / \pi$, $b = B / \pi$, and $c = C / \pi$ to obtain values $a$, $b$, and $c$, where $a + b + c = 1$. We may then define the set of all normalized angle triples by
\begin{equation*}
\overline{\Delta} = \{(a, b, c) \in \mathbb{R}^{3} \mid a, b, c \geq 0, a + b + c = 1\}
\end{equation*}
and the set of \emph{non-degenerate} normalized angle triples by
\begin{equation*}
\Delta = \{(a, b, c) \in \mathbb{R}^{3} \mid a, b, c > 0, a + b + c = 1\};
\end{equation*}
that is, $\Delta$ is the interior of the region defined by $\overline{\Delta}$.

Note that, since we are only concerned with triangle similarity, we regard two triangles having the same angles but in a different order as belonging to the same similarity class. This convention follows that of V\'{a}lyi~\cite{Valyi1903Fusspunktdreiecke}, who considered two triangles to be similar up to a relabelling of the vertices. Thus, to fix an ordering, we will sort each triple of normalized angles in descending order. Similar to before, we define the set of sorted triples by
\begin{equation*}
\overline{C} = \{(a, b, c) \in \mathbb{R}^{3} \mid a \geq b \geq c \geq 0, a + b + c = 1\}
\end{equation*}
and the set of \emph{non-degenerate} sorted triples by
\begin{equation*}
C = \{(a, b, c) \in \mathbb{R}^{3} \mid a \geq b \geq c > 0, a + b + c = 1\}.
\end{equation*}
Any element $(a, b, c) \in C$ therefore refers to the similarity class of triangles whose interior angles are $a\pi$, $b\pi$, and $c\pi$, in any order.


\section{Enumerating Words and Triangles}\label{sec:enumerating}

There exists a formula to enumerate all primitive two-dimensional words of dimension $m \times n$ over a $k$-ary alphabet $\Sigma$:
\begin{equation}
\psi_{k}(m,n) = \sum_{d_{1} \divides m} \sum_{d_{2} \divides n} \mu(d_{1}) \mu(d_{2}) k^{mn / (d_{1}d_{2})}.
\label{eq:2Dprimitiveenum}
\end{equation}
Here, $\mu$ is the M\"{o}bius function~\cite{Mobius1832UberEineBesondere}:
\begin{equation*}
\mu(n) = 
\begin{cases}
1			& \text{if } n = 1; \\
(-1)^{k}	& \text{if } n \text{ is the product of } k \text{ distinct primes; and} \\
0			& \text{if } n \text{ is divisible by a square} > 1.
\end{cases}
\end{equation*}
Equation~\ref{eq:2Dprimitiveenum} is due to Gamard, Richomme, Shallit, and the author~\cite{GamardRichommeShallitSmith2017Periodicity, Smith2017MMathThesis}, and it is essentially a generalization of the well-known formula for enumerating primitive (one-dimensional) words presented in, for example, the book by Lothaire~\cite[p.\ 9]{Lothaire1983CombinatoricsOnWords}.

For the sake of fixing the simplest possible nontrivial case, suppose $\Sigma = \{0, 1\}$ (since, over a unary alphabet, the only primitive two-dimensional word is the unique word of dimension $1 \times 1$) and suppose one of the dimensions is greater than one; say, $m = 2$ (so as to exclude one-dimensional words). Then, we may consider all two-dimensional words consisting of two rows and some number of columns $n \geq 1$. Let us consider a few small illustrative examples.

\begin{example}\label{ex:primitive2Dwords}
The two primitive two-dimensional words of dimension $2 \times 1$ over $\Sigma = \{0, 1\}$ are
\begin{equation*}
\begin{bsmallmatrix}
0 \\ 1
\end{bsmallmatrix}
\text{ and }
\begin{bsmallmatrix}
1 \\ 0
\end{bsmallmatrix}.
\end{equation*}
The ten primitive two-dimensional words of dimension $2 \times 2$ over $\Sigma$ are
\begin{align*}
&\begin{bsmallmatrix}
0 & 0 \\ 0 & 1
\end{bsmallmatrix}, \
\begin{bsmallmatrix}
0 & 0 \\ 1 & 0
\end{bsmallmatrix}, \
\begin{bsmallmatrix}
0 & 1 \\ 0 & 0
\end{bsmallmatrix}, \
\begin{bsmallmatrix}
1 & 0 \\ 0 & 0
\end{bsmallmatrix}, \\
&\begin{bsmallmatrix}
1 & 1 \\ 1 & 0
\end{bsmallmatrix}, \
\begin{bsmallmatrix}
1 & 1 \\ 0 & 1
\end{bsmallmatrix}, \
\begin{bsmallmatrix}
1 & 0 \\ 1 & 1
\end{bsmallmatrix}, \
\begin{bsmallmatrix}
0 & 1 \\ 1 & 1
\end{bsmallmatrix}, \\
&\begin{bsmallmatrix}
0 & 1 \\ 1 & 0
\end{bsmallmatrix},
\text{ and }
\begin{bsmallmatrix}
1 & 0 \\ 0 & 1
\end{bsmallmatrix}.
\end{align*}
\end{example}

\begin{table}[t]
\centering
\caption{Values of $\psi_{2}(m,n)$, the number of primitive two-dimensional words of dimension $m \times n$ over a binary alphabet}
\label{tab:primitive2D}
\begin{tabular}{c|ccccccc}
$m \setminus n$	& 1		& 2			& 3			& 4				& 5					& 6					& $\cdots$ \\
\hline
1				& 2		& 2			& 6			& 12			& 30				& 54				& \\
2				& 2		& 10		& 54		& 228			& 990				& 3\,966			& \\
3				& 6		& 54		& 498		& 4\,020		& 32\,730			& 261\,522			& \\
4				& 12	& 228		& 4\,020	& 65\,040		& 1\,047\,540		& 16\,768\,860		& $\cdots$ \\
5				& 30	& 990		& 32\,730	& 1\,047\,540	& 33\,554\,370		& 1\,073\,708\,010	& \\
6				& 54	& 3\,966	& 261\,522	& 16\,768\,860	& 1\,073\,708\,010	& 68\,718\,945\,018	& \\
$\vdots$		& 		& 			&			& $\vdots$		&					& 					&
\end{tabular}
\end{table}

Table~\ref{tab:primitive2D} lists values of $\psi_{2}(m,n)$ for $1 \leq m \leq 6$ and $1 \leq n \leq 6$. Due to space considerations, the table is constrained to include only the first six values of $m$ and $n$. However, one may find additional values under various sequence numbers in the On-line Encyclopedia of Integer Sequences (OEIS)~\cite{OEIS}: the first row of the table is \seqnum{A027375}, the third row is \seqnum{A395475}, the fourth row is \seqnum{A290754}, the fifth row is \seqnum{A291070}, the sixth row is \seqnum{A291071}, and the diagonal is \seqnum{A265627}.

There also exists a formula to enumerate all triangles having exact pedal period $n$, due to V\'{a}lyi~\cite{Valyi1903Fusspunktdreiecke}. If we take $\phi(n)$ to be the number of all triangles similar to their $n$th pedal triangle,\footnote{V\'{a}lyi uses the notation $\psi(n)$, but we use $\phi$ here to avoid confusion with $\psi_{k}(m,n)$.} then V\'{a}lyi observes that $\phi(n) = 2^{n}(2^{n} - 1)$. However, V\'{a}lyi remarks that this number includes those triangles that are similar to their $d$th pedal triangle, where $d$ is a divisor of $n$. By using inclusion-exclusion to restrict ourselves to counting only triangles whose \emph{first} similar triangle is their $n$th pedal triangle, we obtain the following formula for the number of triangles having exact pedal period $n$:
\begin{equation}
\chi(n) = \phi(n) - \sum_{p_{1} \in P} \phi\left(\frac{n}{p_{1}}\right) + \sum_{p_{1}, p_{2} \in P} \phi\left(\frac{n}{p_{1}p_{2}}\right) - \cdots,
\end{equation}
where $P = \{p \text{ prime} \mid p \divides n\}$ is the set of prime factors of $n$ and each $p_{i}$ in a sum is distinct and unordered.

Again, let us consider a few small illustrative examples.

\begin{example}\label{ex:pedaltriangles}
The two triangles having exact pedal period $1$ (that is, similar to their first pedal triangle) are
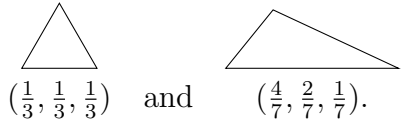
\begin{center}
\begin{tabular}{c c c}
\begin{tikzpicture}
  \def\A{1}
  \def\angA{60} 
  \def\angB{60} 
  \pgfmathsetmacro{\angC}{180-\angA-\angB}
  \pgfmathsetmacro{\d}{\A/sin(\angA)}
  \pgfmathsetmacro{\C}{\d*sin(\angC)}
  \draw (0,0) -- (\angB:\A) -- (0:\C) -- cycle;
\end{tikzpicture} & &
\begin{tikzpicture}
  \def\A{1}
  \def\angA{25} 
  \def\angB{51} 
  \pgfmathsetmacro{\angC}{180-\angA-\angB}
  \pgfmathsetmacro{\d}{\A/sin(\angA)}
  \pgfmathsetmacro{\C}{\d*sin(\angC)}
  \draw (0,0) -- (\angB:\A) -- (0:\C) -- cycle;
\end{tikzpicture} \\
$(\frac{1}{3},\frac{1}{3},\frac{1}{3})$ & and &
$(\frac{4}{7},\frac{2}{7},\frac{1}{7})$.
\end{tabular}
\end{center}
The ten triangles having exact pedal period $2$ (that is, similar to their second pedal triangle, but not to their first pedal triangle) are
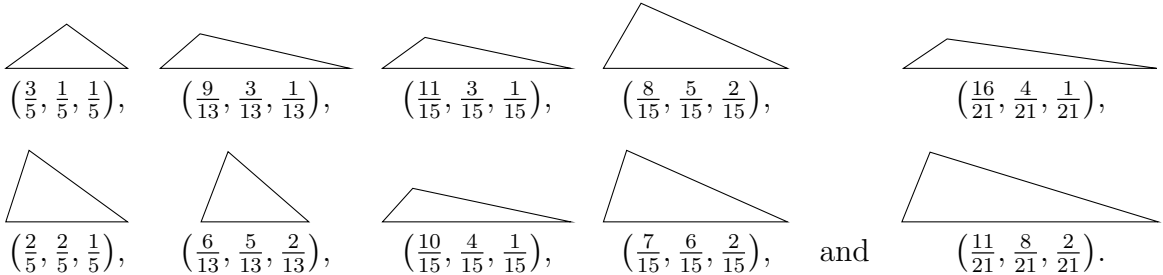
\begin{center}
\begin{tabular}{c c c c c c}
\begin{tikzpicture}
  \def\A{1}
  \def\angA{36} 
  \def\angB{36} 
  \pgfmathsetmacro{\angC}{180-\angA-\angB}
  \pgfmathsetmacro{\d}{\A/sin(\angA)}
  \pgfmathsetmacro{\C}{\d*sin(\angC)}
  \draw (0,0) -- (\angB:\A) -- (0:\C) -- cycle;
\end{tikzpicture} &
\begin{tikzpicture}
  \def\A{0.7}
  \def\angA{13} 
  \def\angB{41} 
  \pgfmathsetmacro{\angC}{180-\angA-\angB}
  \pgfmathsetmacro{\d}{\A/sin(\angA)}
  \pgfmathsetmacro{\C}{\d*sin(\angC)}
  \draw (0,0) -- (\angB:\A) -- (0:\C) -- cycle;
\end{tikzpicture} &
\begin{tikzpicture}
  \def\A{0.7}
  \def\angA{12} 
  \def\angB{36} 
  \pgfmathsetmacro{\angC}{180-\angA-\angB}
  \pgfmathsetmacro{\d}{\A/sin(\angA)}
  \pgfmathsetmacro{\C}{\d*sin(\angC)}
  \draw (0,0) -- (\angB:\A) -- (0:\C) -- cycle;
\end{tikzpicture} &
\begin{tikzpicture}
  \def\A{1}
  \def\angA{24} 
  \def\angB{60} 
  \pgfmathsetmacro{\angC}{180-\angA-\angB}
  \pgfmathsetmacro{\d}{\A/sin(\angA)}
  \pgfmathsetmacro{\C}{\d*sin(\angC)}
  \draw (0,0) -- (\angB:\A) -- (0:\C) -- cycle;
\end{tikzpicture} & &
\begin{tikzpicture}
  \def\A{0.7}
  \def\angA{8} 
  \def\angB{34} 
  \pgfmathsetmacro{\angC}{180-\angA-\angB}
  \pgfmathsetmacro{\d}{\A/sin(\angA)}
  \pgfmathsetmacro{\C}{\d*sin(\angC)}
  \draw (0,0) -- (\angB:\A) -- (0:\C) -- cycle;
\end{tikzpicture} \\
$\left(\frac{3}{5},\frac{1}{5},\frac{1}{5}\right)$, &
$\left(\frac{9}{13},\frac{3}{13},\frac{1}{13}\right)$, &
$\left(\frac{11}{15},\frac{3}{15},\frac{1}{15}\right)$, &
$\left(\frac{8}{15},\frac{5}{15},\frac{2}{15}\right)$, & &
$\left(\frac{16}{21},\frac{4}{21},\frac{1}{21}\right)$, \\[1em]
\begin{tikzpicture}
  \def\A{1}
  \def\angA{36} 
  \def\angB{72} 
  \pgfmathsetmacro{\angC}{180-\angA-\angB}
  \pgfmathsetmacro{\d}{\A/sin(\angA)}
  \pgfmathsetmacro{\C}{\d*sin(\angC)}
  \draw (0,0) -- (\angB:\A) -- (0:\C) -- cycle;
\end{tikzpicture} &
\begin{tikzpicture}
  \def\A{1}
  \def\angA{41} 
  \def\angB{69} 
  \pgfmathsetmacro{\angC}{180-\angA-\angB}
  \pgfmathsetmacro{\d}{\A/sin(\angA)}
  \pgfmathsetmacro{\C}{\d*sin(\angC)}
  \draw (0,0) -- (\angB:\A) -- (0:\C) -- cycle;
\end{tikzpicture} &
\begin{tikzpicture}
  \def\A{0.6}
  \def\angA{12} 
  \def\angB{48} 
  \pgfmathsetmacro{\angC}{180-\angA-\angB}
  \pgfmathsetmacro{\d}{\A/sin(\angA)}
  \pgfmathsetmacro{\C}{\d*sin(\angC)}
  \draw (0,0) -- (\angB:\A) -- (0:\C) -- cycle;
\end{tikzpicture} &
\begin{tikzpicture}
  \def\A{1}
  \def\angA{24} 
  \def\angB{72} 
  \pgfmathsetmacro{\angC}{180-\angA-\angB}
  \pgfmathsetmacro{\d}{\A/sin(\angA)}
  \pgfmathsetmacro{\C}{\d*sin(\angC)}
  \draw (0,0) -- (\angB:\A) -- (0:\C) -- cycle;
\end{tikzpicture} & &
\begin{tikzpicture}
  \def\A{1}
  \def\angA{17} 
  \def\angB{68} 
  \pgfmathsetmacro{\angC}{180-\angA-\angB}
  \pgfmathsetmacro{\d}{\A/sin(\angA)}
  \pgfmathsetmacro{\C}{\d*sin(\angC)}
  \draw (0,0) -- (\angB:\A) -- (0:\C) -- cycle;
\end{tikzpicture} \\
$\left(\frac{2}{5},\frac{2}{5},\frac{1}{5}\right)$, &
$\left(\frac{6}{13},\frac{5}{13},\frac{2}{13}\right)$, &
$\left(\frac{10}{15},\frac{4}{15},\frac{1}{15}\right)$, &
$\left(\frac{7}{15},\frac{6}{15},\frac{2}{15}\right)$, & and &
$\left(\frac{11}{21},\frac{8}{21},\frac{2}{21}\right)$.
\end{tabular}
\end{center}
\end{example}

\begin{table}[t]
\centering
\caption{Values of $\chi(n)$, the number of triangles having exact pedal period $n$}
\label{tab:pedal}
\begin{tabular}{ccccccccccc}
1 & 2 & 3 & 4 & 5 & 6 & 7 & 8 & 9 & 10 & $\cdots$	\\
\hline
2 & 10 & 54 & 228 & 990 & 3\,966 & 16\,254 & 65\,040 & 261\,576 & 1\,046\,550 & $\cdots$
\end{tabular}
\end{table}

Table~\ref{tab:pedal} lists values of $\chi(n)$ for $1 \leq n \leq 10$. Again, only the first ten values are provided here due to space considerations. In the OEIS, sequence \seqnum{A102536} lists additional values.

One may have observed earlier that no sequence number from the OEIS was provided for the second row of Table~\ref{tab:primitive2D}. This is because the second row coincides termwise with the pedal triangle sequence \seqnum{A102536}. We may establish this via the following theorem.

\begin{theorem}\label{thm:counts}
$\psi_{2}(2,n) = \chi(n)$ for all $n \geq 1$.

\begin{proof}
By substituting $k = 2$ and $m = 2$ into Equation~\ref{eq:2Dprimitiveenum}, we obtain
\begin{align*}
\psi_{2}(2,n)	&= \sum_{d_{1} \divides 2} \sum_{d_{2} \divides n} \mu(d_{1}) \mu(d_{2}) 2^{2n / (d_{1}d_{2})} \\
				&= \sum_{d \divides n} \mu(d) 2^{2n/d} - \sum_{d \divides n} \mu(d) 2^{2n/2d} \\
				&= \sum_{d \divides n} \mu(d) \left(4^{n/d} - 2^{n/d}\right).
\end{align*}

Now, recall V\'{a}lyi observed that $\phi(n) = 2^{n}(2^{n} - 1) = 4^{n} - 2^{n}$. However, since this number includes triangles similar to their $d$th pedal triangle, where $d$ is a divisor of $n$, we may alternatively express this number as
\begin{equation}\label{eq:phin}
\phi(n) = \sum_{d \divides n} \chi(d).
\end{equation}
Then, applying M\"{o}bius inversion to Equation~\ref{eq:phin}, we obtain
\begin{align*}
\chi(n)	&= \sum_{d \divides n} \mu(d) \ \phi\left(\frac{n}{d}\right) \\
		&= \sum_{d \divides n} \mu(d) \left(4^{n/d} - 2^{n/d}\right). \qedhere
\end{align*}
\end{proof}
\end{theorem}


\section{Constructing the Bijection}\label{sec:bijection}

Having established via Theorem~\ref{thm:counts} that the counts of our objects of study are the same, we now turn to constructing the bijection between these two sets of objects. Interestingly, nothing about the structure of either object itself leads to the bijection; rather, as we will show, both sets of objects are related by way of the same finite symbolic dynamical system.


\subsection{Two-Dimensional Words to One-Dimensional Words}\label{subsec:2Dto1D}

Let us begin by focusing on two-dimensional words. Suppose $\Sigma = \{0, 1\}$. For $n \geq 1$, let $\mathcal{W}_{n} \subseteq \Sigma^{2 \times n}$ denote the set of primitive two-dimensional words of dimension $2 \times n$ over $\Sigma$. In what follows, we will need the following substitution from alphabet symbols to the language of columns of words in $\mathcal{W}_{n}$.

\begin{definition}\label{def:2D1Dsubstitution}
Let $\Gamma = \{0, 1, 2, 3\}$. Take $\eta \from \Gamma \to \Sigma^{2 \times 1}$ to be as follows:
\begin{equation*}
\eta(0) = \begin{bsmallmatrix}0 \\ 1\end{bsmallmatrix}, \
\eta(1) = \begin{bsmallmatrix}0 \\ 0\end{bsmallmatrix}, \
\eta(2) = \begin{bsmallmatrix}1 \\ 1\end{bsmallmatrix}, \text{ and }
\eta(3) = \begin{bsmallmatrix}1 \\ 0\end{bsmallmatrix}.
\end{equation*}
Extend $\eta$ from symbols to words of length $n$ over $\Gamma$ in the usual way.
\end{definition}

There are two forbidding conditions we must rule out in order for a word over $\Gamma$ to correspond to some word in $\mathcal{W}_{n}$. First, any word over $\Gamma$ of the form $\{1, 2\}^{n}$ cannot map to a primitive two-dimensional word, as applying $\eta$ would produce a two-dimensional word where the first and second rows are identical. Second, if some word $w$ over $\Gamma$ is periodic, then the two-dimensional word $\eta(w)$ would similarly be periodic. In light of these observations, let $\mathcal{U}_{n} \subseteq \Gamma^{n}$ denote the set of primitive \emph{column words} $w$ of length $n$ over $\Gamma$ where $w \not\in \{1, 2\}^{n}$.

\begin{lemma}\label{lem:2D1Dbijection}
The substitution $\eta$ restricts to a bijection between $\mathcal{U}_{n}$ and $\mathcal{W}_{n}$.

\begin{proof}
Let $W \in \Sigma^{2 \times n}$ for some $n \geq 1$, and let $w \in \Gamma^{n}$ be such that $W = \eta(w)$.

First, suppose $W \not\in \mathcal{W}_{n}$; that is, $W$ is a periodic two-dimensional word. Then because $W$ contains two rows, two cases arise:
\begin{itemize}
\item The two rows of $W$ are equal. Then $W = V^{2 \times 1}$ for some $V \in \Sigma^{1 \times n}$, and hence $w \in \{1, 2\}^{n}$.
\item The two rows of $W$ are not equal. Then $W = V^{1 \times q}$ for some $V \in \Sigma^{2 \times p}$, $p \geq 1$, and $q \geq 2$, and $w = v^{q}$ for some $v \in \Gamma^{p}$ such that $V = \eta(v)$.
\end{itemize}

Next, suppose $w \not\in \mathcal{U}_{n}$. If $w \in \{1, 2\}^{n}$, then every column of $W$ is either $\begin{bsmallmatrix}0 \\ 0\end{bsmallmatrix}$ or $\begin{bsmallmatrix}1 \\ 1\end{bsmallmatrix}$, and $W$ is rowwise periodic; while if $w = v^{q}$ for some $q \geq 2$, then $W$ is columnwise periodic.

Lastly, as the substitution $\eta$ is columnwise bijective, it is therefore bijective. The result follows.
\end{proof}
\end{lemma}


\subsection{Sorted Pedal Map}\label{subsec:pedalmap}

Next, we shift our focus to pedal triangles. Hobson~\cite{Hobson1891PlaneTrigonometry} was the first to provide angle formulas for constructing pedal triangles, while Kingston and Synge~\cite{KingstonSynge1988PedalTriangles} corrected errors in Hobson's work and studied the sequences that arose by iteratively constructing pedal triangles according to the corrected angle formulas. Denoting the original triangle by $T$ and its corresponding triple by $(a, b, c) \in \Delta$, the sequence develops according to the following pedal mapping:
\begin{itemize}
\item If $T$ is acute, then its pedal triangle corresponds to the triple $(1-2a, 1-2b, 1-2c)$.
\item If $T$ is obtuse at angle $A$, then its pedal triangle corresponds to the triple $(2a-1, 2b, 2c)$.
\item If $T$ is obtuse at angle $B$, then its pedal triangle corresponds to the triple $(2a, 2b-1, 2c)$.
\item If $T$ is obtuse at angle $C$, then its pedal triangle corresponds to the triple $(2a, 2b, 2c-1)$.
\end{itemize}
Kingston and Synge~\cite{KingstonSynge1988PedalTriangles} observed that the process of iteratively applying this pedal mapping gives rise to a discontinuous four-to-one mapping of a simplex to itself, where each triple $(a, b, c)$ gives the barycentric coordinates of a point interior to this simplex. Lax~\cite{Lax1990ErgodicCharacterPedalTriangles} later showed that this pedal mapping is ergodic, while Ungar~\cite{Ungar1990MixingPropertyPedal} represented this pedal mapping as a shift.

Here, we will use a similar idea to develop a pedal map for our sorted triples in $C$. Before we proceed, let $C^{*} = \{(a, b, c) \in C \mid a \neq 1/2\}$; we may interpret this set $C^{*}$ as the set of sorted triples for which at least one step of the pedal map is defined. We exclude sorted triples where $a = 1/2$ because this only occurs when $T$ is a right triangle and is therefore pedally degenerate.

\begin{definition}\label{def:pedalregions}
Partition the set $C^{*}$ into the following four regions:
\begin{align*}
R_{0}	&= \{(a, b, c) \in C^{*} \mid a < 1/2\}; \\
R_{1}	&= \{(a, b, c) \in C^{*} \mid 2a-1 \geq 2b, a > 1/2\}; \\
R_{2}	&= \{(a, b, c) \in C^{*} \mid 2b > 2a-1 \geq 2c, a > 1/2\}; \text{ and} \\
R_{3}	&= \{(a, b, c) \in C^{*} \mid 2c > 2a-1, a > 1/2\}.
\end{align*}
\end{definition}

\begin{definition}\label{def:pedalmap}
Take $P \from C^{*} \to C$ to be the following sorted pedal map:
\begin{equation*}
P(a, b, c) = 
\begin{cases}
(1-2c, 1-2b, 1-2a)	& \text{if } (a, b, c) \in R_{0}; \\
(2a-1, 2b, 2c)		& \text{if } (a, b, c) \in R_{1}; \\
(2b, 2a-1, 2c)		& \text{if } (a, b, c) \in R_{2}; \text{ and} \\
(2b, 2c, 2a-1)		& \text{if } (a, b, c) \in R_{3}.
\end{cases}
\end{equation*}
\end{definition}


\subsection{Inverse Branches of the Sorted Pedal Map}\label{subsec:inversebranch}

Observe that the map $P$ from Definition~\ref{def:pedalmap} is expanding on each branch. As the remainder of our proof will use a finite symbolic coding, it will be more convenient for us to use inverse branches, which results in a contracting map.

\begin{definition}\label{def:inversebranch}
For each branch of the sorted pedal map $P$, take the corresponding inverse branch on $\overline{C}$ to be as follows:
\begin{align*}
I_{0}(a, b, c)	&= \left( \frac{1-c}{2}, \frac{1-b}{2}, \frac{1-a}{2} \right), \\
I_{1}(a, b, c)	&= \left( \frac{a+1}{2}, \frac{b}{2}, \frac{c}{2} \right), \\
I_{2}(a, b, c)	&= \left( \frac{b+1}{2}, \frac{a}{2}, \frac{c}{2} \right), \text{ and} \\
I_{3}(a, b, c)	&= \left( \frac{c+1}{2}, \frac{a}{2}, \frac{b}{2} \right).
\end{align*}
\end{definition}

Observe that we may compose branches in the following way. For any column word $w = w_{0}w_{1} \dots w_{m-1}$ over $\Gamma$, define $I_{w} = I_{w_{0}} \circ I_{w_{1}} \circ \dots \circ I_{w_{m-1}}$ with composition applied from right to left.

We now establish a few fundamental properties of these inverse branches.

\begin{lemma}\label{lem:inverseproperties}
For each $i \in \{0, 1, 2, 3\}$, the inverse branch $I_{i}$ maps $\overline{C}$ to $\overline{C}$, is injective, and is a contraction with Lipschitz constant $1/2$.

\begin{proof}
To see that an inverse branch $I_{i}$ maps $\overline{C}$ to $\overline{C}$, consider applying $I_{0}$ to some triple $\overline{p} = (a, b, c)$. Since $a \geq b \geq c$, it is the case that $(1-c)/2 \geq (1-b)/2 \geq (1-a)/2$. Additionally,
\begin{equation*}
\frac{1-c}{2} + \frac{1-b}{2} + \frac{1-a}{2} = \frac{3 - (a+b+c)}{2} = 1,
\end{equation*}
and so $I_{0}(\overline{p}) \in \overline{C}$ for all $\overline{p} \in \overline{C}$. By similar arguments, we see that each of $I_{1}$, $I_{2}$, and $I_{3}$ also map $\overline{C}$ to $\overline{C}$.

Each inverse branch $I_{i}$ is injective by observation.

Finally, for any $\overline{p} = (a, b, c)$ and $\overline{p}^{\prime} = (a^{\prime}, b^{\prime}, c^{\prime})$ in $\overline{C}$, a direct calculation shows that $d(I_{i}(\overline{p}), I_{i}(\overline{p}^{\prime}))^{2} = 1/4 \cdot d(\overline{p}, \overline{p}^{\prime})^{2}$, where $d$ is the usual Euclidean distance metric and $i \in \{0, 1, 2, 3\}$. Therefore, $d(I_{i}(\overline{p}), I_{i}(\overline{p}^{\prime})) = 1/2 \cdot d(\overline{p}, \overline{p}^{\prime})$, and so each inverse branch $I_{i}$ is a contraction with Lipschitz constant $1/2$.
\end{proof}
\end{lemma}

The next result reveals why we have temporarily been using $\overline{C}$ in this subsection instead of $C$ by itself: combining $\overline{C}$ with the Euclidean distance metric $d$ gives us a \emph{complete} metric space, which allows us to conclude the following.

\begin{lemma}\label{lem:fixedpoint}
For all $w \in \Gamma^{m}$ where $m \geq 1$, $I_{w} \from \overline{C} \to \overline{C}$ has a unique fixed point.

\begin{proof}
By its definition, $\overline{C}$ is a closed subset of $\mathbb{R}^{3}$, and so it together with the Euclidean distance metric $d$ forms a nonempty, complete metric space. By Lemma~\ref{lem:inverseproperties}, $I_{w}$ is a contraction map with constant $1/2^{m}$, as it is the composition of $m$ contractions each with Lipschitz constant $1/2$. Since $m \geq 1$, it is the case that $1/2^{m} < 1$. Thus, by the Banach fixed-point theorem, there exists a unique point $\overline{p}_{w} \in \overline{C}$ such that $I_{w}(\overline{p}_{w}) = \overline{p}_{w}$.
\end{proof}
\end{lemma}


\subsection{Branch Itineraries}\label{subsec:branchitinerary}

For $n \geq 1$, let $\mathcal{T}_{n}$ denote the set of sorted triples corresponding to triangles having exact pedal period $n$. Using our pedal map notation, we may write this definition as
\begin{equation*}
\mathcal{T}_{n} = \{p \in C \mid P^{j}(p) \text{ is defined for all } 0 \leq j < n, P^{n}(p) = p, \text{ and } P^{d}(p) \neq p \text{ for any } 1 \leq d < n\}.
\end{equation*}

As one of the final steps toward our bijection, we define a so-called \emph{itinerary map} that serves to record the sequence of branches of our sorted pedal map produced by some sorted triple $p$ corresponding to a triangle in $\mathcal{T}_{n}$.

\begin{definition}\label{def:itinerarymap}
Take $\iota \from \mathcal{T}_{n} \to \Gamma^{n}$ to be an itinerary map, and for any $p \in \mathcal{T}_{n}$, take the branch itinerary to be the word $\iota(p) = w_{0}w_{1} \dots w_{n-1}$, where $P^{j}(p) \in R_{w_{j}}$ for all $0 \leq j < n$.
\end{definition}

In other terms, the itinerary map $\iota$ encodes branches arising from the partition given in Definition~\ref{def:pedalregions} as words over a four-symbol alphabet. Specifically, the symbol $w_{j}$ records which of the four branches of the sorted pedal map is used at the $j$th step of the period. This particular idea is not entirely novel; Alexander~\cite{Alexander1993SymbolicDynamicsPedalTriangles} similarly used a four-symbol alphabet to assign to ordered triangles an infinite ``obtuseness sequence", and pedal iteration corresponds to a left shift being applied to this sequence. By contrast, here we sort the triple of angles after each pedal iteration and record via the branch itinerary which branch of the sorted pedal map was used. Thus, we essentially make use of a finite sorted analogue of Alexander's infinite coding.

Before we continue, we require a few technical results pertaining to the itinerary map. The first technical result establishes that, given a branch itinerary $w$, the triple associated with that branch itinerary is the fixed point of the composition of inverse branches $I_{w}$.

\begin{lemma}\label{lem:itinerarymapfixedpoint}
If $p \in \mathcal{T}_{n}$ and $\iota(p) = w$, then $p = I_{w}(p)$.

\begin{proof}
Suppose, as before, that $w = w_{0}w_{1} \dots w_{n-1}$. Since $P^{j}(p) \in R_{w_{j}}$ for all $0 \leq j < n$, the inverse branch $I_{w_{j}}$ allows us to recover $P^{j}(p)$ from $P^{j+1}(p)$. Therefore, we can write
\begin{equation*}
p = I_{w_{0}}(P(p)) = I_{w_{0}} \circ I_{w_{1}}(P^{2}(p)) = \dots = I_{w_{0}} \circ I_{w_{1}} \circ \dots \circ I_{w_{n-1}}(P^{n}(p)) = I_{w}(P^{n}(p)).
\end{equation*}
Since $p \in \mathcal{T}_{n}$, we have that $P^{n}(p) = p$, and so $p = I_{w}(p)$.
\end{proof}
\end{lemma}

The second technical result connects branch itineraries of sorted triples to words over $\Gamma$ of the form $\{1, 2\}^{n}$; recall that this is one of the forbidding conditions we discussed in Section~\ref{subsec:2Dto1D}.

\begin{lemma}\label{lem:itinerarymapprimitive1}
If $p \in \mathcal{T}_{n}$, then $\iota(p) \not\in \{1, 2\}^{n}$.

\begin{proof}
Let $p = (a, b, c) \in \mathcal{T}_{n}$ and suppose by way of contradiction that $\iota(p) \in \{1, 2\}^{n}$. Then at each point, we must follow either the branch corresponding to $R_{1}$ or the branch corresponding to $R_{2}$. On both of these branches, the third coordinate is doubled: $P(a, b, c) = (2a-1, 2b, 2c)$ in the case of $R_{1}$, or $P(a, b, c) = (2b, 2a-1, 2c)$ in the case of $R_{2}$. After $n$ iterations, the third coordinate becomes $2^{n}c$. However, since $p \in \mathcal{T}_{n}$, we know that $P^{n}(p) = p$. Hence, $c = 2^{n}c$. Since $n \geq 1$, this forces $c = 0$, which means that $p \not\in C$ and therefore $p \not\in \mathcal{T}_{n}$.
\end{proof}
\end{lemma}

The third technical result connects branch itineraries of sorted triples to the notion of periodicity, which again is one of the forbidding conditions mentioned in Section~\ref{subsec:2Dto1D}.

\begin{lemma}\label{lem:itinerarymapprimitive2}
If $p \in \mathcal{T}_{n}$, then $\iota(p)$ is a primitive word over $\Gamma$.

\begin{proof}
Let $w = \iota(p)$, and suppose by way of contradiction that $w$ is periodic. Then $w = v^{m}$ for some $v \in \Gamma^{\ell}$ with $\ell \geq 1$ and $m \geq 2$. Since $w = \iota(p)$, we have that $n = \ell m$ where $\ell < n$.

By Lemma~\ref{lem:itinerarymapfixedpoint}, we know that $p = I_{w}(p)$, and since $w = v^{m}$, we may rewrite the inverse branch as $I_{w} = (I_{v})^{m}$ and therefore $p = (I_{v})^{m}(p)$. By Lemma~\ref{lem:fixedpoint}, we know that $I_{v}$ has a unique fixed point. Additionally, since $I_{v}$ is contractive, $(I_{v})^{m}$ is also contractive and, by the Banach fixed-point theorem, has the same fixed point as $I_{v}$. Thus, we may write $p = I_{v}(p)$ by itself.

Since the first $\ell$ symbols of the branch itinerary are $v$ itself, performing the inverse branch calculation demonstrated in the proof of Lemma~\ref{lem:itinerarymapfixedpoint} will produce $p = I_{v}(P^{\ell}(p))$. By Lemma~\ref{lem:inverseproperties}, we know that $I_{v}$ is injective, and so taking both $p = I_{v}(p)$ and $p = I_{v}(P^{\ell}(p))$ together implies that $p = P^{\ell}(p)$. However, this means that $p$ cannot have exact pedal period $n$, and therefore $p \not\in \mathcal{T}_{n}$.
\end{proof}
\end{lemma}

Lastly, we combine the previous results to establish some helpful properties about the itinerary map.

\begin{lemma}\label{lem:itinerarymapproperties}
The itinerary map $\iota \from \mathcal{T}_{n} \to \mathcal{U}_{n}$ is well defined and injective.

\begin{proof}
Let $p \in \mathcal{T}_{n}$. By Definition~\ref{def:itinerarymap}, we have that $\iota(p) \in \Gamma^{n}$. Since the four regions $R_{0}$, $R_{1}$, $R_{2}$, and $R_{3}$ form a disjoint partition of the domain of the sorted pedal map according to Definition~\ref{def:pedalregions}, each iteration of the pedal map $P^{j}(p)$ falls into a unique region $R_{i}$, and hence $\iota(p)$ is uniquely determined. Then, by Lemma~\ref{lem:itinerarymapprimitive1} and Lemma~\ref{lem:itinerarymapprimitive2}, we know both that $\iota(p) \not\in \{1, 2\}^{n}$ and that $\iota(p)$ is a primitive word, so $\iota(p) \in \mathcal{U}_{n}$. Therefore, $\iota$ is well defined as a map from $\mathcal{T}_{n}$ to $\mathcal{U}_{n}$.

Now, consider $p \in \mathcal{T}_{n}$ and $p^{\prime} \in \mathcal{T}_{n}$ such that $\iota(p) = \iota(p^{\prime}) = w$. By Lemma~\ref{lem:itinerarymapfixedpoint}, we know that $p = I_{w}(p)$ and $p^{\prime} = I_{w}(p^{\prime})$, and so both $p$ and $p^{\prime}$ are fixed points of $I_{w}$. However, Lemma~\ref{lem:fixedpoint} tells us that $I_{w}$ must have a unique fixed point in $\overline{C}$, and so it must be the case that $p = p^{\prime}$. Thus, $\iota$ is injective.
\end{proof}
\end{lemma}


\subsection{Main Theorem}\label{subsec:bijectiontheorem}

We now have enough to prove the main theorem of the paper.

\begin{theorem}\label{thm:bijection}
For each $n \geq 1$, the map $\eta \circ \iota$ is a bijection between $\mathcal{T}_{n}$ and $\mathcal{W}_{n}$.

\begin{proof}
By Lemma~\ref{lem:itinerarymapproperties}, we know that $\iota \from \mathcal{T}_{n} \to \mathcal{U}_{n}$ is injective. Furthermore, we know by Lemma~\ref{lem:2D1Dbijection} that the substitution $\eta$ restricts to a bijection between $\mathcal{U}_{n}$ and $\mathcal{W}_{n}$. Therefore, by composition, we have that $\eta \circ \iota \from \mathcal{T}_{n} \to \mathcal{W}_{n}$ is injective.

Now, by Theorem~\ref{thm:counts}, we know that $\psi_{2}(2,n) = \chi(n)$ for all $n \geq 1$. We also know that $|\mathcal{W}_{n}| = \psi_{2}(2,n)$ and $|\mathcal{T}_{n}| = \chi(n)$ by definition, so the sets $\mathcal{W}_{n}$ and $\mathcal{T}_{n}$ have the same cardinality.

Since $\eta \circ \iota$ is an injection between two finite sets having the same cardinality, it is a bijection as desired.
\end{proof}
\end{theorem}

It is a rather straightforward exercise to illustrate, using Theorem~\ref{thm:bijection}, that a given triangle having exact pedal period $n$ determines a unique primitive two-dimensional word of dimension $2 \times n$ over a binary alphabet, and that such a two-dimensional word itself corresponds to a unique triangle having exact pedal period $n$. Starting with a triangle having exact pedal period $n$, for instance, we sort the triple corresponding to the triangle to obtain $p \in \mathcal{T}_{n}$, record the branch itinerary, and then apply the substitution $\eta$. Conversely, starting with a two-dimensional word $W \in \mathcal{W}_{n}$, we compute $w = \eta^{-1}(W)$, then find the associated triangle obtained by solving the fixed-point equation $I_{w}(p) = p$ in $\overline{C}$.

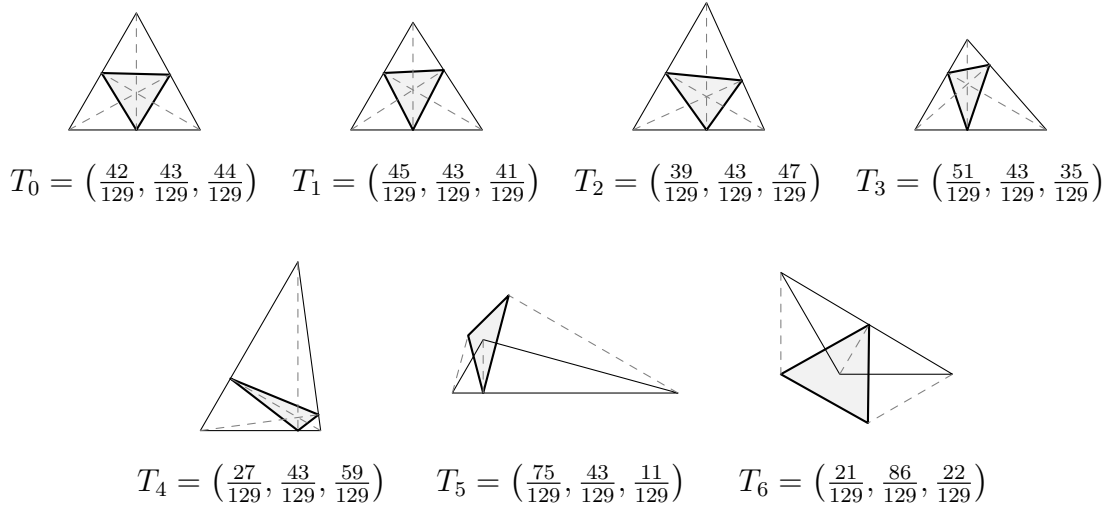
\begin{figure}[t]
\centering
\begin{tikzpicture}[line cap=round,line join=round]
\newcommand{\pedal}[9]{
	\begin{scope}[shift={#1}]
		\pgfmathsetmacro{\sum}{int(round(#3+#4+#5))}
		\pgfmathsetmacro{\edge}{#9}
		\pgfmathsetmacro{\Adeg}{180*(#3)/\sum}
		\pgfmathsetmacro{\Bdeg}{180*(#4)/\sum}
		\pgfmathsetmacro{\Cdeg}{180*(#5)/\sum}

		\pgfmathsetmacro{\AB}{\edge*sin(\Cdeg)/sin(\Adeg)}
		\pgfmathsetmacro{\Ax}{\AB*cos(\Bdeg)}
		\pgfmathsetmacro{\Ay}{\AB*sin(\Bdeg)}
		\pgfmathsetmacro{\Bx}{0}
		\pgfmathsetmacro{\By}{0}
		\pgfmathsetmacro{\Cx}{\edge}
		\pgfmathsetmacro{\Cy}{0}

		\pgfmathsetmacro{\Xx}{\Ax}
		\pgfmathsetmacro{\Xy}{0}

		\pgfmathsetmacro{\dYx}{\Ax-\Cx}
		\pgfmathsetmacro{\dYy}{\Ay-\Cy}
		\pgfmathsetmacro{\tY}{((\Bx-\Cx)*\dYx+(\By-\Cy)*\dYy)/(\dYx*\dYx+\dYy*\dYy)}
		\pgfmathsetmacro{\Yx}{\Cx+\tY*\dYx}
		\pgfmathsetmacro{\Yy}{\Cy+\tY*\dYy}

		\pgfmathsetmacro{\dZx}{\Bx-\Ax}
		\pgfmathsetmacro{\dZy}{\By-\Ay}
		\pgfmathsetmacro{\tZ}{((\Cx-\Ax)*\dZx+(\Cy-\Ay)*\dZy)/(\dZx*\dZx+\dZy*\dZy)}
		\pgfmathsetmacro{\Zx}{\Ax+\tZ*\dZx}
		\pgfmathsetmacro{\Zy}{\Ay+\tZ*\dZy}

		\coordinate (A) at (\Ax,\Ay);
		\coordinate (B) at (\Bx,\By);
		\coordinate (C) at (\Cx,\Cy);
		\coordinate (X) at (\Xx,\Xy);
		\coordinate (Y) at (\Yx,\Yy);
		\coordinate (Z) at (\Zx,\Zy);

		\draw[black, thick, fill=gray!10] (X)--(Y)--(Z)--cycle;
		\draw (A)--(B)--(C)--cycle;

		\draw[gray, dashed] (A)--(X);
		\draw[gray, dashed] (B)--(Y);
		\draw[gray, dashed] (C)--(Z);

		\node[align=center] at #8 {$#2=\left(\frac{#3}{\sum},\frac{#4}{\sum},\frac{#5}{\sum}\right)$};
	\end{scope}
}
\pedal{(0,0)}{T_0}{42}{43}{44}{T_1}{(45,43,41)}{(0.88,-0.7)}{1.75}
\pedal{(3.75,0)}{T_1}{45}{43}{41}{T_2}{(39,43,47)}{(0.88,-0.7)}{1.75}
\pedal{(7.5,0)}{T_2}{39}{43}{47}{T_3}{(51,43,35)}{(0.88,-0.7)}{1.75}
\pedal{(11.25,0)}{T_3}{51}{43}{35}{T_4}{(27,43,59)}{(0.88,-0.7)}{1.75}
\pedal{(1.75,-4)}{T_4}{27}{43}{59}{T_5}{(75,43,11)}{(0.82,-0.7)}{1.6}
\pedal{(5.1,-3.5)}{T_5}{75}{43}{11}{T_6}{(21,86,22)}{(1.45,-1.2)}{3}
\pedal{(10.25,-3.25)}{T_6}{21}{86}{22}{T_7}{(42,43,44)}{(0.32,-1.45)}{1.5}
\end{tikzpicture}
\caption{A pedal heptacycle given by Kingston and Synge~\cite{KingstonSynge1988PedalTriangles}. Each pedal triangle $T_{i+1}$ (in grey) is constructed from the triangle $T_{i}$ (in white) for $0 \leq i \leq 6$. The triples corresponding to each triangle are listed underneath the triangles.}
\label{fig:kingstonsynge}
\end{figure}

As an example, consider the following pedal heptacycle given by Kingston and Synge~\cite{KingstonSynge1988PedalTriangles}:
\begin{align*}
&\left(\frac{42}{129}, \frac{43}{129}, \frac{44}{129}\right) \rightarrow
\left(\frac{45}{129}, \frac{43}{129}, \frac{41}{129}\right) \rightarrow
\left(\frac{39}{129}, \frac{43}{129}, \frac{47}{129}\right) \rightarrow
\left(\frac{51}{129}, \frac{43}{129}, \frac{35}{129}\right) \\ \rightarrow \
&\left(\frac{27}{129}, \frac{43}{129}, \frac{59}{129}\right) \rightarrow
\left(\frac{75}{129}, \frac{43}{129}, \frac{11}{129}\right) \rightarrow
\left(\frac{21}{129}, \frac{86}{129}, \frac{22}{129}\right) \rightarrow
\left(\frac{42}{129}, \frac{43}{129}, \frac{44}{129}\right).
\end{align*}
While an equilateral triangle has exact pedal period $1$, as demonstrated in Example~\ref{ex:pedaltriangles}, this nearly equilateral triangle has exact pedal period $7$: by iterating the pedal map seven times on the initial sorted triple $\left(\frac{44}{129}, \frac{43}{129}, \frac{42}{129}\right)$, we obtain a similar triangle represented by the same triple. Figure~\ref{fig:kingstonsynge} illustrates the heptacycle as given by Kingston and Synge; observe that the pedal triangle $T_{7}$ produced from the triangle $T_{6}$ is similar to the original triangle $T_{0}$. In what follows, we will use the same triples, but in their sorted forms.

\begin{example}
Let $T_{0}$ be the triangle corresponding to the triple $\left(\frac{42}{129}, \frac{43}{129}, \frac{44}{129}\right)$. Sort this triple in descending order to obtain $p = \left(\frac{44}{129}, \frac{43}{129}, \frac{42}{129}\right)$.

Since $44/129 < 1/2$, this triple will be in region $R_{0}$. The first branch of the pedal map gives us
\begin{align*}
P\left(\frac{44}{129}, \frac{43}{129}, \frac{42}{129}\right) &= \left(1-2\left(\frac{42}{129}\right), 1-2\left(\frac{43}{129}\right), 1-2\left(\frac{44}{129}\right)\right) \\
&= \left(\frac{45}{129}, \frac{43}{129}, \frac{41}{129}\right),
\end{align*}
as expected. Indeed, one can verify that each of the first five triples $p$, $P(p)$, $P^{2}(p)$, $P^{3}(p)$, and $P^{4}(p)$ place us in region $R_{0}$.

Now, observe that $P^{5}(p) = \left(\frac{75}{129}, \frac{43}{129}, \frac{11}{129}\right)$, and since $75/129 > 1/2$, we must be in one of the other three regions. Since $2c = 22/129 > 21/129 = 2a-1$, this triple will be in region $R_{3}$.

Applying the pedal map again produces $P^{6}(p) = \left(\frac{86}{129}, \frac{22}{129}, \frac{21}{129}\right)$; we have that $86/129 > 1/2$, and since $2b = 44/129 > 43/129 = 2a-1$ and $2a-1 = 43/129 \geq 42/129 = 2c$, this triple will be in region $R_{2}$. Applying the pedal map once more produces $P^{7}(p) = \left(\frac{44}{129}, \frac{43}{129}, \frac{42}{129}\right) = p$.

From these iterated applications of the pedal map, we find that the branch itinerary is $\iota(p) = 0000032$. We now apply the substitution $\eta$ to obtain
\begin{align*}
\eta(\iota(p))	&= \eta(0000032) \\
				&= \eta(0) \, \eta(0) \, \eta(0) \, \eta(0) \, \eta(0) \, \eta(3) \, \eta(2) \\
				&= \begin{bsmallmatrix}0 & 0 & 0 & 0 & 0 & 1 & 1 \\ 1& 1 & 1 & 1 & 1 & 0 & 1\end{bsmallmatrix},
\end{align*}
and this is the primitive two-dimensional word $W$ associated with the triangle $T_{0}$ having exact pedal period $7$.
\end{example}


\section*{Acknowledgements}\label{sec:acknowledgements}

This research was supported by the Natural Sciences and Engineering Research Council of Canada (NSERC) Discovery Grant RGPIN-2024-04799.


\bibliographystyle{plain}
\bibliography{References}


\end{document}